\documentclass{amsart}
\usepackage{amssymb}
\usepackage{epsfig}

\title{Schur positivity and Schur log-concavity}

\author{Thomas Lam, Alexander Postnikov, Pavlo Pylyavskyy}

\date{February 18, 2005, updated September 26, 2005}

\address{Department of Mathematics, M.I.T., Cambridge, MA 02139}

\email{thomasl (at) math (dot) mit (dot) edu}
\email{apost (at) math (dot) mit (dot) edu}
\email{pasha (at) math (dot) mit (dot) edu}

\keywords{Schur functions, Schur positivity, Schur log-concavity, immanants,
Kazhdan-Lusztig polynomials, Temperley-Lieb algebra, minors}

\thanks{A.P.\ was supported in part by NSF grant DMS-0201494.  P.P.\ would like to thank Institut Mittag-Leffler for their hospitality.}

\theoremstyle{plain}
\newtheorem{theorem}{Theorem}

\newtheorem{corollary}[theorem]{Corollary}
\newtheorem{conjecture}[theorem]{Conjecture}

\theoremstyle{definition}

\theoremstyle{remark}

\newcommand{\la}[1]{\overleftarrow{#1}}
\newcommand{\ra}[1]{\overrightarrow{#1}}

\def\ainv{\mathrm{ainv}}

\def\uqsln{U_q(\widehat{ \mathfrak{
sl}}_n)}
\def\ll{\lambda}

\def\C{\mathbb C}

\def\X{\,\,\lower2pt\hbox{
\begin{picture}(0,0)%
\includegraphics{figX.pstex}%
\end{picture}%
\setlength{\unitlength}{1973sp}%
\begingroup\makeatletter\ifx\SetFigFont\undefined%
\gdef\SetFigFont#1#2#3#4#5{%
  \reset@font\fontsize{#1}{#2pt}%
  \fontfamily{#3}\fontseries{#4}\fontshape{#5}%
  \selectfont}%
\fi\endgroup%
\begin{picture}(324,324)(589,-973)
\end{picture}
}}
\def\noXv{\,\,\lower2pt\hbox{
\begin{picture}(0,0)%
\includegraphics{figNoX.pstex}%
\end{picture}%
\setlength{\unitlength}{1973sp}%
\begingroup\makeatletter\ifx\SetFigFont\undefined%
\gdef\SetFigFont#1#2#3#4#5{%
  \reset@font\fontsize{#1}{#2pt}%
  \fontfamily{#3}\fontseries{#4}\fontshape{#5}%
  \selectfont}%
\fi\endgroup%
\begin{picture}(316,316)(293,-969)
\end{picture}
}}
\def\noXh{\,\,\lower2pt\hbox{
\begin{picture}(0,0)%
\includegraphics{figNoX2.pstex}%
\end{picture}%
\setlength{\unitlength}{1973sp}%
\begingroup\makeatletter\ifx\SetFigFont\undefined%
\gdef\SetFigFont#1#2#3#4#5{%
  \reset@font\fontsize{#1}{#2pt}%
  \fontfamily{#3}\fontseries{#4}\fontshape{#5}%
  \selectfont}%
\fi\endgroup%
\begin{picture}(316,316)(893,-369)
\end{picture}
}}
\def\noXDU{\,\,\lower2pt\hbox{
\begin{picture}(0,0)%
\includegraphics{figDU.pstex}%
\end{picture}%
\setlength{\unitlength}{1973sp}%
\begingroup\makeatletter\ifx\SetFigFont\undefined%
\gdef\SetFigFont#1#2#3#4#5{%
  \reset@font\fontsize{#1}{#2pt}%
  \fontfamily{#3}\fontseries{#4}\fontshape{#5}%
  \selectfont}%
\fi\endgroup%
\begin{picture}(316,316)(1493,-969)
\end{picture}
}}
\def\noXDD{\,\,\lower2pt\hbox{
\begin{picture}(0,0)%
\includegraphics{figDD.pstex}%
\end{picture}%
\setlength{\unitlength}{1973sp}%
\begingroup\makeatletter\ifx\SetFigFont\undefined%
\gdef\SetFigFont#1#2#3#4#5{%
  \reset@font\fontsize{#1}{#2pt}%
  \fontfamily{#3}\fontseries{#4}\fontshape{#5}%
  \selectfont}%
\fi\endgroup%
\begin{picture}(316,316)(1493,-969)
\end{picture}
}}
\def\noXUD{\,\,\lower2pt\hbox{
\begin{picture}(0,0)%
\includegraphics{figUD.pstex}%
\end{picture}%
\setlength{\unitlength}{1973sp}%
\begingroup\makeatletter\ifx\SetFigFont\undefined%
\gdef\SetFigFont#1#2#3#4#5{%
  \reset@font\fontsize{#1}{#2pt}%
  \fontfamily{#3}\fontseries{#4}\fontshape{#5}%
  \selectfont}%
\fi\endgroup%
\begin{picture}(316,316)(1493,-969)
\end{picture}
}}
\def\noXUU{\,\,\lower2pt\hbox{
\begin{picture}(0,0)%
\includegraphics{figUU.pstex}%
\end{picture}%
\setlength{\unitlength}{1973sp}%
\begingroup\makeatletter\ifx\SetFigFont\undefined%
\gdef\SetFigFont#1#2#3#4#5{%
  \reset@font\fontsize{#1}{#2pt}%
  \fontfamily{#3}\fontseries{#4}\fontshape{#5}%
  \selectfont}%
\fi\endgroup%
\begin{picture}(316,316)(1493,-969)
\end{picture}
}}
\def\noXRR{\,\,\lower2pt\hbox{
\begin{picture}(0,0)%
\includegraphics{figRR.pstex}%
\end{picture}%
\setlength{\unitlength}{1973sp}%
\begingroup\makeatletter\ifx\SetFigFont\undefined%
\gdef\SetFigFont#1#2#3#4#5{%
  \reset@font\fontsize{#1}{#2pt}%
  \fontfamily{#3}\fontseries{#4}\fontshape{#5}%
  \selectfont}%
\fi\endgroup%
\begin{picture}(316,316)(1493,-969)
\end{picture}
}}
\def\noXRL{\,\,\lower2pt\hbox{
\begin{picture}(0,0)%
\includegraphics{figRL.pstex}%
\end{picture}%
\setlength{\unitlength}{1973sp}%
\begingroup\makeatletter\ifx\SetFigFont\undefined%
\gdef\SetFigFont#1#2#3#4#5{%
  \reset@font\fontsize{#1}{#2pt}%
  \fontfamily{#3}\fontseries{#4}\fontshape{#5}%
  \selectfont}%
\fi\endgroup%
\begin{picture}(316,316)(1493,-969)
\end{picture}
}}
\def\noXLR{\,\,\lower2pt\hbox{
\begin{picture}(0,0)%
\includegraphics{figLR.pstex}%
\end{picture}%
\setlength{\unitlength}{1973sp}%
\begingroup\makeatletter\ifx\SetFigFont\undefined%
\gdef\SetFigFont#1#2#3#4#5{%
  \reset@font\fontsize{#1}{#2pt}%
  \fontfamily{#3}\fontseries{#4}\fontshape{#5}%
  \selectfont}%
\fi\endgroup%
\begin{picture}(316,316)(1493,-969)
\end{picture}
}}
\def\noXLL{\,\,\lower2pt\hbox{
\begin{picture}(0,0)%
\includegraphics{figLL.pstex}%
\end{picture}%
\setlength{\unitlength}{1973sp}%
\begingroup\makeatletter\ifx\SetFigFont\undefined%
\gdef\SetFigFont#1#2#3#4#5{%
  \reset@font\fontsize{#1}{#2pt}%
  \fontfamily{#3}\fontseries{#4}\fontshape{#5}%
  \selectfont}%
\fi\endgroup%
\begin{picture}(316,316)(1493,-969)
\end{picture}
}}
\def\XRR{\,\,\lower2pt\hbox{
\begin{picture}(0,0)%
\includegraphics{figXRR.pstex}%
\end{picture}%
\setlength{\unitlength}{1973sp}%
\begingroup\makeatletter\ifx\SetFigFont\undefined%
\gdef\SetFigFont#1#2#3#4#5{%
  \reset@font\fontsize{#1}{#2pt}%
  \fontfamily{#3}\fontseries{#4}\fontshape{#5}%
  \selectfont}%
\fi\endgroup%
\begin{picture}(324,324)(589,-973)
\end{picture}%
}}
\def\XRL{\,\,\lower2pt\hbox{
\begin{picture}(0,0)%
\includegraphics{figXRL.pstex}%
\end{picture}%
\setlength{\unitlength}{1973sp}%
\begingroup\makeatletter\ifx\SetFigFont\undefined%
\gdef\SetFigFont#1#2#3#4#5{%
  \reset@font\fontsize{#1}{#2pt}%
  \fontfamily{#3}\fontseries{#4}\fontshape{#5}%
  \selectfont}%
\fi\endgroup%
\begin{picture}(324,324)(589,-973)
\end{picture}%
}}
\def\XLR{\,\,\lower2pt\hbox{
\begin{picture}(0,0)%
\includegraphics{figXLR.pstex}%
\end{picture}%
\setlength{\unitlength}{1973sp}%
\begingroup\makeatletter\ifx\SetFigFont\undefined%
\gdef\SetFigFont#1#2#3#4#5{%
  \reset@font\fontsize{#1}{#2pt}%
  \fontfamily{#3}\fontseries{#4}\fontshape{#5}%
  \selectfont}%
\fi\endgroup%
\begin{picture}(324,324)(589,-973)
\end{picture}%
}}
\def\XLL{\,\,\lower2pt\hbox{
\begin{picture}(0,0)%
\includegraphics{figXLL.pstex}%
\end{picture}%
\setlength{\unitlength}{1973sp}%
\begingroup\makeatletter\ifx\SetFigFont\undefined%
\gdef\SetFigFont#1#2#3#4#5{%
  \reset@font\fontsize{#1}{#2pt}%
  \fontfamily{#3}\fontseries{#4}\fontshape{#5}%
  \selectfont}%
\fi\endgroup%
\begin{picture}(324,324)(589,-973)
\end{picture}%
}}

\def\Imm{\mathrm{Imm}}
\def\TL{\mathit{TL}}

\begin{document}

\begin{abstract}
We prove Okounkov's conjecture, a conjecture of Fomin-Fulton-Li-Poon,
and a special case of Lascoux-Leclerc-Thibon's conjecture on
Schur positivity and give several more general statements
using a recent result of Rhoades and Skandera. 
An alternative proof of this result is provided.
We also give an intriguing log-concavity property of Schur functions.
\end{abstract}


\maketitle

\section{Schur positivity conjectures}
\label{sec:snn}

The ring of symmetric functions has a linear basis of {\it Schur
functions} $s_\lambda$ labelled by partitions $\lambda= (\lambda_1\geq
\lambda_2\geq \dots\geq 0)$, see~\cite{Mac}.  These functions appear in
representation theory as characters of irreducible representations of $GL_n$
and in geometry as representatives of Schubert classes for Grassmannians.  A
symmetric function is called {\it Schur nonnegative} if it is a linear
combination with nonnegative coefficients of the Schur functions,
or, equivalently,  
if it is the character of a certain
representation of $GL_n$.  In particular, {\it skew Schur functions}
$s_{\lambda/\mu}$ are Schur nonnegative.  
Recently, a lot of work has gone into studying whether
certain expressions of the form $s_\ll s_\mu - s_\nu s_\rho$ were Schur
nonnegative.  Schur positivity of an expression of this form is equivalent
to some inequalities between Littlewood-Richardson coefficients.
In a sense, characterizing such inequalities 
is a ``higher analogue'' of the Klyachko problem on nonzero
Littlewood-Richardson coefficients.
Let us 
mention several Schur positivity conjectures due to Okounkov,
Fomin-Fulton-Li-Poon, and Lascoux-Leclerc-Thibon of the above form.

Okounkov~\cite{Oko} studied branching rules for classical Lie
groups and proved that the multiplicities were ``monomial log-concave''
in some sense.
An essential combinatorial ingredient in his construction was
the theorem about monomial nonnegativity of some symmetric functions.
He conjectured that these functions are Schur nonnegative, as well.
For a partition $\ll$ with all even parts, let $\frac{\ll}{2}$ denote
the partition $(\frac{\ll_1}{2},\frac{\ll_2}{2},\ldots)$.
For two symmetric functions $f$ and $g$, the notation $f\geq_s g$ means that $f-g$ is
Schur nonnegative.  

\begin{conjecture}
\label{conj:Oko}
{\rm Okounkov~\cite{Oko}} \
For two skew shapes $\ll/\mu$ and  $\nu/\rho$ such that
$\ll+\nu$ and $\mu+\rho$
both have all even parts,  we have
$(s_{\frac{(\ll + \nu)}{2}/\frac{(\mu+\rho)}{2}})^2 \geq_s
s_{\ll/\mu}\, s_{\nu/\rho}$.
\end{conjecture}

Fomin, Fulton, Li, and Poon \cite{FFLP} studied the eigenvalues
and singular values of sums of Hermitian and of complex matrices.
Their study led to two combinatorial conjectures concerning
differences of products of Schur functions. Let us formulate one
of these conjectures, which was also studied recently by Bergeron
and McNamara~\cite{BM}.  For two partitions $\ll$ and $\mu$, let
$\ll \cup \mu=(\nu_1,\nu_2,\nu_3,\dots)$ be the partition obtained
by rearranging all parts of $\lambda$ and $\mu$ in the weakly
decreasing order. Let $\mathrm{sort}_1(\lambda,\mu) :=
(\nu_1,\nu_3,\nu_5,\dots)$ and $\mathrm{sort}_2(\lambda,\mu) :=
(\nu_2,\nu_4,\nu_6,\dots)$.

\begin{conjecture}
\label{conj:FFLP}
{\rm Fomin-Fulton-Li-Poon \cite[Conjecture~2.7]{FFLP}} \
For two partitions $\ll$ and $\mu$, we have
$s_{\mathrm{sort}_1(\ll,\mu)} s_{\mathrm{sort}_2(\ll,\mu)} \geq_s s_\ll s_\mu$.
\end{conjecture}

Lascoux, Leclerc, and Thibon \cite{LLT} studied a family of
symmetric functions $\mathcal G^{(n)}_{\ll}(q,x)$
arising combinatorially from ribbon tableaux
and algebraically from the Fock space representations of the
quantum affine algebra $\uqsln$.
They conjectured that ${\mathcal G}^{(n)}_{n\ll}(q,x) \geq_s
{\mathcal G}^{(m)}_{m\ll}(q,x)$ for $m \leq n$.
For the case $q = 1$,
their conjecture can be reformulated, as follows.
For a partition $\lambda$ and $1\leq i\leq n$, let
$\lambda^{[i,n]}:=(\lambda_i,\lambda_{i+n},\lambda_{i+2n},\dots)$.
In particular, $\mathrm{sort}_i(\lambda,\mu) = (\lambda\cup\mu)^{[i,2]}$,
for $i=1,2$.

\begin{conjecture}
\label{conj:LLT}
{\rm Lascoux-Leclerc-Thibon~\cite[Conjecture~6.4]{LLT}} \
For integers $1 \leq m \leq n$ and a partition $\ll$, we have
$\prod_{i=1}^n
s_{\ll^{[i,n]}} \geq_s \prod_{i=1}^m s_{\ll^{[i,m]}}$.
\end{conjecture}

\begin{theorem}  Conjectures~\ref{conj:Oko}, \ref{conj:FFLP} and~\ref{conj:LLT}
are true.
\end{theorem}

In Section~\ref{sec:proofs}, we present and prove more general versions
of these conjectures.
Our approach is based on the following result.
For two partitions $\lambda=(\lambda_1, \lambda_2, \dots)$ and
$\mu=(\mu_1, \mu_2, \dots)$, let us define partitions
$\lambda\vee \mu:=(\max(\lambda_1,\mu_1),\max(\lambda_2,\mu_2),\dots)$
and
$\lambda\wedge \mu:=(\min(\lambda_1,\mu_1),\min(\lambda_2,\mu_2),\dots)$.
The Young diagram of $\lambda\vee\mu$ is the set-theoretical union
of the Young diagrams of $\lambda$ and $\mu$.
Similarly, the Young diagram of $\lambda\wedge\mu$ is the set-theoretical
intersection of the Young diagrams of $\lambda$ and $\mu$.
For two skew shapes, define
$(\lambda/\mu)\vee(\nu/\rho) :=
\lambda\vee\nu/\mu\vee\rho$  and
$(\lambda/\mu)\wedge(\nu/\rho) :=
\lambda\wedge\nu/\mu\wedge\rho$.

\begin{theorem}
\label{th:nonnegativity_s}
Let $\lambda/\mu$ and $\nu/\rho$ be any two skew shapes.
Then we have
$$
s_{(\lambda/\mu) \vee (\nu/\rho)}\,
s_{(\lambda/\mu) \wedge (\nu/\rho)} \geq_s
s_{\lambda/\mu}\, s_{\nu/\rho}.
$$
\end{theorem}

This theorem was originally conjectured by Lam and Pylyavskyy in~\cite{LP}.

\section{Background}

In this section we give an overview of some results 
of Haiman~\cite{Hai} and Rhoades-Skandera~\cite{RS, RS2}.
We include an alternative proof Rhoades-Skandera's result.

\subsection{Haiman's Schur positivity result}

Let $H_n(q)$ be the {\it Hecke algebra\/} associated with the
symmetric group $S_n$.
The Hecke algebra has the standard basis $\{T_w \mid w \in S_n\}$
and the {\it Kazhdan-Lusztig basis\/} $\{C_w'(q)\mid w\in S_n\}$ related by
\[
q^{l(v)/2} C_v'(q) = \sum_{w \leq v} P_{w,v}(q) \,T_w
\quad\textrm{and}\quad T_w = \sum_{v \leq w} (-1)^{l(vw)} Q_{v,w}(q)
\, q^{l(v)/2} C'_v(q),
\]
where $P_{w,v}(q)$ are the {\it Kazhdan-Lusztig polynomials\/} and
$Q_{v,w}(q) = P_{w_\circ w,w_\circ v}(q)$, for the longest
permutation $w_\circ\in S_n$, see~\cite{Hum} for more details. 

For $w \in S_n$ and a $n \times n$ matrix $X = (x_{ij})$,
the {\it Kazhdan-Lusztig immanant\/} was defined in~\cite{RS}
as
$$
\Imm_{w}(X):=\sum_{v \in S_n} (-1)^{l(vw)}Q_{w,v}(1)
\,x_{1,v(1)} \cdots x_{n,v(n)},
$$  

Let $h_k = \sum_{i_1\leq \cdots \leq i_k} x_{i_1}\cdots x_{i_k}$ 
be the $k$-th
homogeneous symmetric function, where $h_0 = 1$ and $h_k = 0$ for $k
< 0$.  A {\it generalized Jacobi-Trudi matrix\/} is a $n\times n$ matrix of the form 
$\left(h_{\mu_i-\nu_j}\right)_{i,j = 1}^n$,
for partitions $\mu = (\mu_1 \geq \mu_2 \cdots \geq \mu_n \geq 0)$ and
$\nu = (\nu_1 \geq \nu_2 \cdots \geq \nu_n \geq 0)$.
Haiman's result
can be reformulated as follows, see~\cite{RS}.

\begin{theorem}
\label{thm:Haiman}
{\rm Haiman~\cite[Theorem~1.5]{Hai}} \ 
The immanants $\Imm_w$ of
a generalized Jacobi-Trudi matrix are Schur non-negative.
\end{theorem}

Haiman's proof of this result is based on the Kazhdan-Lusztig conjecture 
proven by Beilinson-Bernstein and Brylinski-Kashiwara.
This conjecture expresses the characters of 
Verma modules as sums of the characters of some irreducible
highest weight representations of $\mathfrak{sl}_n$ with multiplicities 
equal to $P_{w,v}(1)$.  One can derive from this conjecture 
that the coefficients of Schur functions in $\Imm_w$ are certain tensor 
product multiplicities of irreducible representations.

\subsection{Temperley-Lieb algebra}

The {\it {Temperley-Lieb algebra}} $\TL_n(\xi)$ is the $\mathbb
C[\xi]$-algebra generated by $t_1, \ldots, t_{n-1}$
subject to the relations $t_i^2=\xi\, t_i$, $t_i t_j t_i=t_i$ if
$|i-j|=1$, $t_i t_j = t_j t_i$ if $|i-j| \geq 2$. The dimension of
$\TL_n(\xi)$ equals the $n$-th Catalan number
$C_n=\frac{1}{n+1}\binom{2n}{n}$. 
A {\it 321-avoiding permutation\/} is a permutation $w\in S_n$ that 
has no reduced decomposition of the form $w = \cdots s_i s_j s_i \cdots$ with $|i-j|=1$.
(These permutations are also called {\it fully-commutative.})
A natural basis of the Temperley-Lieb algebra is 
$\{t_w\mid w \textrm{ is a 321-avoiding permutation in } S_n\}$,
where $t_w := t_{i_1} \cdots t_{i_l}$, for a reduced decomposition $w = s_{i_1}
\cdots s_{i_l}$.

The map $\theta:T_{s_i} \mapsto t_i-1$ determines
a homomorphism $\theta:H_n(1)=\C[S_n]\to \TL_n(2)$.
Indeed, the elements $t_i-1$ in $\TL_n(2)$ satisfy the Coxeter relations.

\begin{theorem} {\rm Fan-Green~\cite{FG}} \
\label{th:FG}
The homomorphism $\theta$ acts on the Kazhdan-Lusztig basis $\{C_w'(1)\}$ 
of $H_n(1)$ as follows:
$$
\theta(C_w'(1)) =
\begin{cases}
t_w & \text{if $w$ is $321$-avoiding,}\\
0 & \text{otherwise.}
\end{cases}
$$
\end{theorem}


For any permutation $v \in S_n$ and a 321-avoiding permutation $w\in S_n$, let $f_{w}(v)$
be the coefficient of the basis element $t_w\in \TL_n(2)$ in the
basis expansion of $\theta(T_v) = (t_{i_1}-1)\cdots
(t_{i_l}-1)\in \TL_n(2)$, for a reduced decomposition $v=s_{i_1}\cdots s_{i_l}$.
Rhoades and Skandera~\cite{RS2} defined the {\it Temperley-Lieb
immanant\/} $\Imm_{w}^{\mathrm{TL}}(x)$ of an $n\times n$ matrix $X = (x_{ij})$ by
$$
\Imm_{w}^{\mathrm{TL}}(X):=\sum_{v \in S_n} f_{w}(v)\,x_{1,v(1)} \cdots
x_{n,v(n)}.
$$

\begin{theorem}
\label{th:immschur} 
{\rm Rhoades-Skandera~\cite{RS2}} \
For a 321-avoiding permutation $w\in S_n$, we have
$\Imm_w^{\mathrm{TL}}(X) = \Imm_w(X)$.
\end{theorem}

\begin{proof}
Applying the map $\theta$ to $T_v = \sum_{w \leq v} (-1)^{l(vw)}
Q_{w,v}(1) \, C'_w(1)$ and using Theorem~\ref{th:FG} we obtain
$\theta(T_v)  = \sum
(-1)^{l(vw)} Q_{w,v}(1) \, t_w$,
where the sum is over 321-avoiding permutations $w$.
Thus $f_w(v) = (-1)^{l(vw)} Q_{w,v}(1)$ and
$\Imm_{w}^{\mathrm{TL}}=\Imm_{w}$. 
\end{proof}

A product of generators (decomposition) 
$t_{i_1}\cdots t_{i_l}$ in the Temperley-Lieb algebra $\TL_n$ can be
graphically presented by a {\it Temperley-Lieb diagram\/} with $n$ non-crossing
strands connecting the vertices $1,\dots,2n$ and, possibly, with
some internal loops. This diagram is obtained from the wiring diagram of the
decomposition $w=s_{i_1}\cdots s_{i_l}\in S_n$ by replacing each crossing
``\X'' with a 
{\it vertical uncrossing\/}
``\noXv''.  For example, the following
figure shows the wiring diagram for $s_1 s_2 s_2 s_3 s_2\in S_4$ and the
Temperley-Lieb diagram for $t_1 t_2 t_2 t_3 t_2\in \TL_4$.  
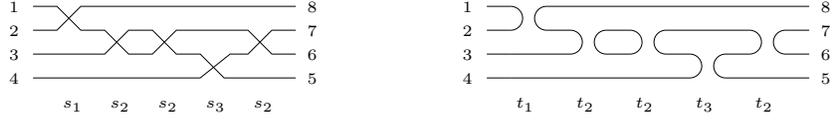
\smallskip

\begin{center}
\begin{picture}(0,0)%
\includegraphics{fig12.pstex}%
\end{picture}%
\setlength{\unitlength}{1973sp}%
\begingroup\makeatletter\ifx\SetFigFont\undefined%
\gdef\SetFigFont#1#2#3#4#5{%
  \reset@font\fontsize{#1}{#2pt}%
  \fontfamily{#3}\fontseries{#4}\fontshape{#5}%
  \selectfont}%
\fi\endgroup%
\begin{picture}(10476,1390)(-5099,-4094)
\put(-4424,-4036){\makebox(0,0)[lb]{\smash{{\SetFigFont{6}{7.2}{\familydefault}{\mddefault}{\updefault}{\color[rgb]{0,0,0}$s_1$}%
}}}}
\put(601,-3736){\makebox(0,0)[lb]{\smash{{\SetFigFont{6}{7.2}{\familydefault}{\mddefault}{\updefault}{\color[rgb]{0,0,0}$4$}%
}}}}
\put(601,-3436){\makebox(0,0)[lb]{\smash{{\SetFigFont{6}{7.2}{\familydefault}{\mddefault}{\updefault}{\color[rgb]{0,0,0}$3$}%
}}}}
\put(601,-3136){\makebox(0,0)[lb]{\smash{{\SetFigFont{6}{7.2}{\familydefault}{\mddefault}{\updefault}{\color[rgb]{0,0,0}$2$}%
}}}}
\put(601,-2836){\makebox(0,0)[lb]{\smash{{\SetFigFont{6}{7.2}{\familydefault}{\mddefault}{\updefault}{\color[rgb]{0,0,0}$1$}%
}}}}
\put(5101,-3736){\makebox(0,0)[lb]{\smash{{\SetFigFont{6}{7.2}{\familydefault}{\mddefault}{\updefault}{\color[rgb]{0,0,0}$5$}%
}}}}
\put(5101,-3436){\makebox(0,0)[lb]{\smash{{\SetFigFont{6}{7.2}{\familydefault}{\mddefault}{\updefault}{\color[rgb]{0,0,0}$6$}%
}}}}
\put(5101,-3136){\makebox(0,0)[lb]{\smash{{\SetFigFont{6}{7.2}{\familydefault}{\mddefault}{\updefault}{\color[rgb]{0,0,0}$7$}%
}}}}
\put(5101,-2836){\makebox(0,0)[lb]{\smash{{\SetFigFont{6}{7.2}{\familydefault}{\mddefault}{\updefault}{\color[rgb]{0,0,0}$8$}%
}}}}
\put(4276,-4036){\makebox(0,0)[lb]{\smash{{\SetFigFont{6}{7.2}{\familydefault}{\mddefault}{\updefault}{\color[rgb]{0,0,0}$t_2$}%
}}}}
\put(3526,-4036){\makebox(0,0)[lb]{\smash{{\SetFigFont{6}{7.2}{\familydefault}{\mddefault}{\updefault}{\color[rgb]{0,0,0}$t_3$}%
}}}}
\put(2776,-4036){\makebox(0,0)[lb]{\smash{{\SetFigFont{6}{7.2}{\familydefault}{\mddefault}{\updefault}{\color[rgb]{0,0,0}$t_2$}%
}}}}
\put(2026,-4036){\makebox(0,0)[lb]{\smash{{\SetFigFont{6}{7.2}{\familydefault}{\mddefault}{\updefault}{\color[rgb]{0,0,0}$t_2$}%
}}}}
\put(1276,-4036){\makebox(0,0)[lb]{\smash{{\SetFigFont{6}{7.2}{\familydefault}{\mddefault}{\updefault}{\color[rgb]{0,0,0}$t_1$}%
}}}}
\put(-5099,-2836){\makebox(0,0)[lb]{\smash{{\SetFigFont{6}{7.2}{\familydefault}{\mddefault}{\updefault}{\color[rgb]{0,0,0}$1$}%
}}}}
\put(-5099,-3136){\makebox(0,0)[lb]{\smash{{\SetFigFont{6}{7.2}{\familydefault}{\mddefault}{\updefault}{\color[rgb]{0,0,0}$2$}%
}}}}
\put(-5099,-3436){\makebox(0,0)[lb]{\smash{{\SetFigFont{6}{7.2}{\familydefault}{\mddefault}{\updefault}{\color[rgb]{0,0,0}$3$}%
}}}}
\put(-5099,-3736){\makebox(0,0)[lb]{\smash{{\SetFigFont{6}{7.2}{\familydefault}{\mddefault}{\updefault}{\color[rgb]{0,0,0}$4$}%
}}}}
\put(-1349,-2836){\makebox(0,0)[lb]{\smash{{\SetFigFont{6}{7.2}{\familydefault}{\mddefault}{\updefault}{\color[rgb]{0,0,0}$8$}%
}}}}
\put(-1349,-3136){\makebox(0,0)[lb]{\smash{{\SetFigFont{6}{7.2}{\familydefault}{\mddefault}{\updefault}{\color[rgb]{0,0,0}$7$}%
}}}}
\put(-1349,-3436){\makebox(0,0)[lb]{\smash{{\SetFigFont{6}{7.2}{\familydefault}{\mddefault}{\updefault}{\color[rgb]{0,0,0}$6$}%
}}}}
\put(-1349,-3736){\makebox(0,0)[lb]{\smash{{\SetFigFont{6}{7.2}{\familydefault}{\mddefault}{\updefault}{\color[rgb]{0,0,0}$5$}%
}}}}
\put(-2024,-4036){\makebox(0,0)[lb]{\smash{{\SetFigFont{6}{7.2}{\familydefault}{\mddefault}{\updefault}{\color[rgb]{0,0,0}$s_2$}%
}}}}
\put(-2624,-4036){\makebox(0,0)[lb]{\smash{{\SetFigFont{6}{7.2}{\familydefault}{\mddefault}{\updefault}{\color[rgb]{0,0,0}$s_3$}%
}}}}
\put(-3224,-4036){\makebox(0,0)[lb]{\smash{{\SetFigFont{6}{7.2}{\familydefault}{\mddefault}{\updefault}{\color[rgb]{0,0,0}$s_2$}%
}}}}
\put(-3824,-4036){\makebox(0,0)[lb]{\smash{{\SetFigFont{6}{7.2}{\familydefault}{\mddefault}{\updefault}{\color[rgb]{0,0,0}$s_2$}%
}}}}
\end{picture}%

\end{center}

Pairs of vertices connected by strands of a wiring diagram are
$(2n+1-i, w(i))$, for $i=1,\dots,n$.
Pairs of vertices connected by strands in a Temperley-Lieb diagram form
a {\it non-crossing matching}, i.e., a graph on the vertices
$1,\dots,2n$ with $n$ disjoint edges that
contains no pair of edges $(a,c)$ and $(b,d)$ with $a<b<c<d$. 
If two Temperley-Lieb diagrams give the same matching and have 
the same number of internal loops,
then the corresponding products of generators of $\TL_n$
are equal to each other.  If the diagram of $a$ is obtained from the diagram of
$b$ by removing $k$ internal loops, then $b = \xi^k a$ in $\TL_n$.

The map that sends $t_w$ to the non-crossing matching
given by its Temperley-Lieb diagram is a bijection between
basis elements $t_w$ of $\TL_n$,  where $w$ is 321-avoiding,
and non-crossing matchings on the vertex set $[2n]$.
For example, the basis element $t_1 t_3 t_2$ of $\TL_4$ 
corresponds to the non-crosssing 
matching with the edges $(1,2),(3,4),(5,8),(6,7)$.


\subsection{An identity for products of minors}
\label{ssec:identity_minors}

For a subset $S\subset [2n]$, let us say that a Temperley-Lieb
diagram (or the associated element in $\TL_n$) is {\it $S$-compatible\/} 
if each strand of the diagram has one end-point in $S$ and the
other end-point in its complement $[2n]\setminus S$.  
Coloring vertices in $S$ black 
and the remaining vertices
white, a basis element $t_w$ is $S$-compatible if and only if each edge 
in the associated matching has two vertices of different colors.  
Let $\Theta(S)$ denote the set of all 321-avoiding permutation
$w\in S_n$ such that $t_w$ is $S$-compatible.

For two subsets $I, J\subset [n]$ of the same cardinality let $\Delta_{I,J}(X)$
denote the {\it minor\/} of an $n\times n$ matrix $X$ in the row set $I$ and the
column set $J$.
Let $\bar I := [n]\setminus I$ and let $I^\wedge := \{2n+1-i\mid i\in I\}$.


\begin{theorem} 
\label{th:immdecomp}
{\rm Rhoades-Skandera~\cite[Proposition~4.3]{RS2}, cf.~Skandera~\cite{Ska}} \
For two subsets $I,J\subset [n]$ of the same cardinality
and $S=J\cup (\bar I)^\wedge$, we have
$$
\Delta_{I,J}(X)\cdot \Delta_{\bar I, \bar J}(X)
= \sum_{w \in \Theta(S)} \Imm_w^{\mathrm{TL}}(X).
$$
\end{theorem}

The proof given in~\cite{RS2} employs planar networks.  
We give a more direct proof that uses the involution principle.

\begin{proof}
Let us fix a permutation $v\in S_n$ with a reduced decomposition
$v = s_{i_1}\cdots s_{i_l}$.
The coefficient of the monomial $x_{1,v(1)}\cdots x_{n,v(n)}$ 
in the expansion of the product 
of two minors $\Delta_{I,J}(X) \cdot \Delta_{\bar I, \bar J}(X)$ equals
$$
\left\{
\begin{array}{cl}
(-1)^{\mathrm{inv}(I)+\mathrm{inv}(\bar I)}&\textrm{if } v(I) = J,\\
0 & \textrm{if } v(I)\ne J,
\end{array}
\right.
$$
where $\mathrm{inv}(I)$ is the number of inversions $i<j$, $v(i)>v(j)$ 
such that $i,j\in I$.

On the other hand, by the definition of $\Imm_w^{\mathrm{TL}}$, the coefficient of 
$x_{1,v(1)}\cdots x_{n,v(n)}$  in the right-hand side of the 
identity equals the sum $\sum (-1)^{r}\, 2^s$ over all diagrams obtained
from the wiring diagram of the reduced decomposition $s_{i_1}\cdots s_{i_l}$
by replacing each crossing ``\X''
with either a {\it vertical uncrossing\/} ``\noXv'' or
a {\it horizontal uncrossing\/} ``\noXh'' 
so that the resulting diagram is $S$-compatible,
where $r$ is the number of horizontal uncrossings ``\noXh''
and $s$ is the number of internal loops in the resulting diagram.
Indeed, the choice of ``\noXv'' corresponds to the 
choice of ``$t_{i_k}$'' and the choice of ``\noXh'' 
corresponds to the choice of ``$-1$'' in the $k$-th term of the product 
$(t_{i_1}-1)\cdots (t_{i_l}-1)\in \TL_n(2)$,
for $k=1,\dots,l$.

Let us pick directions of all strands and loops in such diagrams so that the initial vertex in
each strand belongs to $S$ (and, thus, the end-point is not in $S$).
There are $2^s$ ways to pick directions of $s$ internal loops. 
Thus the above sum can be written as the sum $\sum (-1)^r$ over such 
{\it directed Temperley-Lieb\/} diagrams.

Here is an example of a directed diagram for $v=s_3 s_2 s_1 s_3 s_2 s_3$ 
and $S=\{1,4,5, 7\}$ corresponding
to the term $t_3 t_2 (-1) t_3 (-1) t_3$ in the expansion of the product 
$(t_3-1)(t_2-1)(t_1-1)(t_3-1)(t_2-1)(t_3-1)$.
This diagram comes with the sign $(-1)^2$.

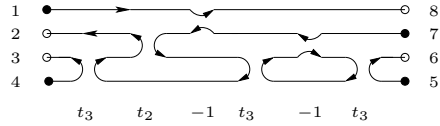
\begin{center}
\begin{picture}(0,0)%
\includegraphics{fig2.pstex}%
\end{picture}%
\setlength{\unitlength}{1973sp}%
\begingroup\makeatletter\ifx\SetFigFont\undefined%
\gdef\SetFigFont#1#2#3#4#5{%
  \reset@font\fontsize{#1}{#2pt}%
  \fontfamily{#3}\fontseries{#4}\fontshape{#5}%
  \selectfont}%
\fi\endgroup%
\begin{picture}(5526,1465)(1,-3269)
\put(3601,-3211){\makebox(0,0)[lb]{\smash{{\SetFigFont{6}{7.2}{\familydefault}{\mddefault}{\updefault}{\color[rgb]{0,0,0}$-1$}%
}}}}
\put(4276,-3211){\makebox(0,0)[lb]{\smash{{\SetFigFont{6}{7.2}{\familydefault}{\mddefault}{\updefault}{\color[rgb]{0,0,0}$t_3$}%
}}}}
\put(5251,-2836){\makebox(0,0)[lb]{\smash{{\SetFigFont{6}{7.2}{\familydefault}{\mddefault}{\updefault}{\color[rgb]{0,0,0}$5$}%
}}}}
\put(5251,-2536){\makebox(0,0)[lb]{\smash{{\SetFigFont{6}{7.2}{\familydefault}{\mddefault}{\updefault}{\color[rgb]{0,0,0}$6$}%
}}}}
\put(5251,-2236){\makebox(0,0)[lb]{\smash{{\SetFigFont{6}{7.2}{\familydefault}{\mddefault}{\updefault}{\color[rgb]{0,0,0}$7$}%
}}}}
\put(5251,-1936){\makebox(0,0)[lb]{\smash{{\SetFigFont{6}{7.2}{\familydefault}{\mddefault}{\updefault}{\color[rgb]{0,0,0}$8$}%
}}}}
\put(1576,-3211){\makebox(0,0)[lb]{\smash{{\SetFigFont{6}{7.2}{\familydefault}{\mddefault}{\updefault}{\color[rgb]{0,0,0}$t_2$}%
}}}}
\put(2251,-3211){\makebox(0,0)[lb]{\smash{{\SetFigFont{6}{7.2}{\familydefault}{\mddefault}{\updefault}{\color[rgb]{0,0,0}$-1$}%
}}}}
\put(2851,-3211){\makebox(0,0)[lb]{\smash{{\SetFigFont{6}{7.2}{\familydefault}{\mddefault}{\updefault}{\color[rgb]{0,0,0}$t_3$}%
}}}}
\put(  1,-1936){\makebox(0,0)[lb]{\smash{{\SetFigFont{6}{7.2}{\familydefault}{\mddefault}{\updefault}{\color[rgb]{0,0,0}$1$}%
}}}}
\put(  1,-2236){\makebox(0,0)[lb]{\smash{{\SetFigFont{6}{7.2}{\familydefault}{\mddefault}{\updefault}{\color[rgb]{0,0,0}$2$}%
}}}}
\put(  1,-2536){\makebox(0,0)[lb]{\smash{{\SetFigFont{6}{7.2}{\familydefault}{\mddefault}{\updefault}{\color[rgb]{0,0,0}$3$}%
}}}}
\put(  1,-2836){\makebox(0,0)[lb]{\smash{{\SetFigFont{6}{7.2}{\familydefault}{\mddefault}{\updefault}{\color[rgb]{0,0,0}$4$}%
}}}}
\put(826,-3211){\makebox(0,0)[lb]{\smash{{\SetFigFont{6}{7.2}{\familydefault}{\mddefault}{\updefault}{\color[rgb]{0,0,0}$t_3$}%
}}}}
\end{picture}%

\end{center}

Let us construct a sign reversing partial involution $\iota$ on the set of 
such directed Temperley-Lieb diagrams.
If a diagram has a {\it misaligned uncrossing\/}, i.e., 
an uncrossing of the form 
``\noXDU'',
``\noXUD'',
``\noXLR'', or
``\noXRL'',
then $\iota$ switches the leftmost such uncrossing according to the rules
$\iota: \textrm{\noXDU}\leftrightarrow \textrm{\noXRL}$
and
$\iota: \textrm{\noXUD}\leftrightarrow \textrm{\noXLR}$.
Otherwise, when the diagram involves only {\it aligned uncrossings\/} 
``\noXUU'',
``\noXDD'',
``\noXLL'',
``\noXRR'',
the involution $\iota$ is not defined.

For example, in the above diagram, the involution $\iota$ switches the second 
uncrossing,
which has the form ``\,\,\lower2pt\hbox{
\begin{picture}(0,0)%
\includegraphics{figUD.pstex}%
\end{picture}%
\setlength{\unitlength}{1973sp}%
\begingroup\makeatletter\ifx\SetFigFont\undefined%
\gdef\SetFigFont#1#2#3#4#5{%
  \reset@font\fontsize{#1}{#2pt}%
  \fontfamily{#3}\fontseries{#4}\fontshape{#5}%
  \selectfont}%
\fi\endgroup%
\begin{picture}(316,316)(1493,-969)
\end{picture}
}'', to
``\,\,\lower2pt\hbox{
\begin{picture}(0,0)%
\includegraphics{figLR.pstex}%
\end{picture}%
\setlength{\unitlength}{1973sp}%
\begingroup\makeatletter\ifx\SetFigFont\undefined%
\gdef\SetFigFont#1#2#3#4#5{%
  \reset@font\fontsize{#1}{#2pt}%
  \fontfamily{#3}\fontseries{#4}\fontshape{#5}%
  \selectfont}%
\fi\endgroup%
\begin{picture}(316,316)(1493,-969)
\end{picture}
}''.  The resulting diagram corresponds
to the term $t_3 (-1) (-1) t_3 (-1) t_3$.

Since the involution $\iota$ reverses signs, this shows that the total 
contribution of all diagrams with at least one misaligned uncrossing is zero.  
Let us show that there is at most one
$S$-compatible directed Temperley-Lieb diagram with all aligned uncrossings.  
If we have a such diagram,
then we can direct the strands of the wiring diagram for $v=s_{i_1}\dots
s_{i_l}$ so that each segment of the wiring diagram has the same direction as
in the Temperley-Lieb diagram.  In particular, the end-points of strands in the
wiring diagram should have different colors.  Thus each strand 
starting at an element of $J$ should finish at an element of $I^\wedge$,
or, equivalently, $v(I)=J$.
The directed Temperley-Lieb diagram can be uniquely recovered from
this directed wiring diagram by replacing the crossings with uncrossings,
as follows: 
$\textrm{\XRR}\to\textrm{\noXRR}$,
$\textrm{\XRL}\to\textrm{\noXUU}$,
$\textrm{\XLR}\to\textrm{\noXDD}$,
$\textrm{\XLL}\to\textrm{\noXLL}$.
Thus the coefficient of 
$x_{1,v(1)}\cdots x_{n,v(n)}$  in the right-hand side
of the needed identity is zero, if $v(I)\ne J$,
and is $(-1)^r$, if $v(I)=J$, where $r$ is the number
of crossings of the form ``\XLL'' or ``\XRR''
in the wiring diagram.
In other words, $r$ equals the number of crossings
such that the right end-points of the pair of crossing strands have 
the same color.
This is exactly the same as the expression for the coefficient
in the left-hand side of the needed identity.
\end{proof}

\section{Proof of Theorem~\ref{th:nonnegativity_s}}
\label{sec:RS}

For two subsets $I,J\subseteq [n]$ of
the same cardinality, let $\Delta_{I,J}(H)$ denote
the minor of the Jacobi-Trudi matrix $H=(h_{j-i})_{1\leq i,j\leq
n}$ with row set $I$ and column set $J$, where $h_i$ is the
$i$-th homogeneous symmetric function, as before.
According to the Jacobi-Trudi formula, see~\cite{Mac}, the minors
$\Delta_{I,J}(H)$ are precisely the skew Schur functions
$$
\Delta_{I,J}(H)= s_{\lambda/\mu},
$$
where $\lambda=(\lambda_1\geq \cdots \geq \lambda_k\geq 0)$,
$\mu=(\mu_1\geq \cdots \geq \mu_k\geq 0)$ and the associated subsets are
$I=\{\mu_k+1<\mu_{k-1}+2<\cdots<\mu_1+k\}$,
$J=\{\lambda_k+1<\lambda_{k-1}+2<\cdots<\lambda_1+k\}$.

For two sets $I=\{i_1<\cdots<i_k\}$ and
$J=\{j_1<\cdots<j_k\}$, let us define
$I\vee J := \{\max(i_1,j_1)< \dots < \max(i_k,j_k)\}$
and $I\wedge J := \{\min(i_1,j_1)< \dots < \min(i_k,j_k)\}$.

%


Theorem~\ref{th:nonnegativity_s} can be reformulated in terms of
minors, as follows.  Without loss of generality we can assume that
all partitions $\lambda,\mu,\nu,\rho$ in
Theorem~\ref{th:nonnegativity_s} have the same number $k$ of
parts, some of which might be zero. 
Note that generalized Jacobi-Trudi matrices are obtained from
$H$ by skipping or duplicating rows and columns.

\begin{theorem}
\label{th:nonnegativity_Delta}
Let $I,J, I', J'$ be $k$ element subsets in $[n]$.
Then we have
$$
\Delta_{I\vee I',\, J\vee J'}(X)\cdot \Delta_{I\wedge I',\, J\wedge J'}(X)  \geq_s
\Delta_{I,J}(X)\cdot \Delta_{I',J'}(X),
$$
for a generalized Jacobi-Trudi matrix $X$.
\end{theorem}

\begin{proof}
Let us denote $\bar I:=[n]\setminus I$ and $\check S:=[2n]\setminus S$.
By skipping or duplicating rows and columns of the matrix $X$, we
may assume that $I'=\bar I$ and $J'=\bar J$. 
Then $I\vee I' = \overline{I\wedge I'}$
and $J\vee J' = \overline{J\wedge J'}$.
Let $S := J \cup  (\bar I)^\wedge$ 
and $T:=(J\vee J') \cup (\overline{ I\vee I'})^\wedge$.
Then we have $T= S \vee \check S$ and $\check T = S\wedge \check S$.

Let us show that $\Theta(S)\subseteq \Theta(T)$, i.e.,
every $S$-compatible non-crossing matching on $[2n]$ 
is also $T$-compatible.
Let $S=\{s_1<\cdots < s_n\}$ and 
$\check S=\{\check s_1 <\cdots < \check s_n\}$.
Then $T=\{\max(s_1,\check s_1),\dots, \max(s_n,\check s_n)\}$ and 
$\check T=\{\min(s_1,\check s_1),\dots, \min(s_n,\check s_n)\}$. 
Let $M$ be an $S$-compatible non-crossing matching on $[2n]$ and let
$(a<b)$ be an edge of $M$.  Without loss of generality we may assume that 
$a = s_i\in S$ and $b = \check s_j \in \check S$.  
We must show that either ($a \in T$ and $b \in \check T$) or 
($a \in \check T$ and $b \in T$).
Since no edge of $M$ can cross $(a,b)$, the elements of $S$ in the 
interval $[a+1,b-1]$ are matched with the elements of $\check S$ 
in this interval. Let $k=\# (S\cap [a+1,b-1]) = \# (\check S\cap [a+1,b-1]$).
Suppose that $a,b\in T$, or, equivalently, $\check s_i < s_i$
and $s_j < \check s_j$.
Since there are at least $k$ elements of $\check S$ in the interval
$[\check s_i + 1, \check s_j -1]$, we have $i+k+1\leq j$.
On the other hand, since there are at most $k-1$ elements of $S$ 
in the interval
$[s_i+1, s_j-1]$, we have $i+k\geq j$.  We obtain a contradiction.
The case $a,b\in\check T$ is analogous.

Now Theorem~\ref{th:immdecomp} implies that 
the difference 
$\Delta_{I\vee I',\, J\vee J'}\cdot \Delta_{I\wedge I',\, J\wedge J'}
- \Delta_{I,J}\cdot \Delta_{I',J'}$ is a nonnegative combination
of Temperley-Lieb immanants.
Theorems~\ref{thm:Haiman} and~\ref{th:immschur} imply
its Schur nonnegativity.
\end{proof}

\section{Proof of conjectures and generalizations}
\label{sec:proofs}

In this section we prove generalized versions of
Conjectures~\ref{conj:Oko}-\ref{conj:LLT}, which were conjectured by
Kirillov~\cite[Section~6.8]{Kir}.  Corollary~\ref{cor:twoshapesconjugate} was
also conjectured by Bergeron-McNamara~\cite[Conjecture~5.2]{BM} who showed that
it implies Theorem~\ref{thm:skewlogconcavity2}.

Let $\lfloor x\rfloor$ be the maximal integer $\leq x$ and $\lceil
x\rceil$ be the minimal integer $\geq x$. For vectors $v$ and $w$
and a positive integer $n$, we assume that the operations $v+w$,
$\frac v n$, $\lfloor v\rfloor$, $\lceil v\rceil$ are performed
coordinate-wise. In particular, we have well-defined operations
$\lfloor\frac{\ll+\nu}{2}\rfloor$ and
$\lceil\frac{\ll+\nu}{2}\rceil$ on pairs of partitions.

The next claim extends Okounkov's conjecture
(Conjecture~\ref{conj:Oko}).

\begin{theorem}
\label{th:twoshapes}
Let $\ll/\mu$ and $\nu/\rho$ be any two skew
shapes. Then we have
\[
s_{\lfloor\frac{\ll+\nu}{2}\rfloor/\lfloor\frac{\mu+\rho}{2}\rfloor}\,
s_{\lceil\frac{\ll+\nu}{2}\rceil/\lceil\frac{\mu+\rho}{2}\rceil}
\geq_s s_{\ll/\mu}\, s_{\nu/\rho}.
\]
\end{theorem}

\begin{proof}
We will assume that all partitions have the same fixed number $k$ of parts,
some of which might be zero.
For a skew shape $\ll/\mu = (\ll_1,\ldots,\ll_k)/(\mu_1,\ldots,\mu_k)$,
define
$$
\ra{\ll/\mu}:= (\ll_1 +1 ,\ldots,\ll_k +1 )/(\mu_1+1,\ldots,\mu_k+1),
$$
that is,  $\ra{\ll/\mu}$ is the skew shape obtained by shifting
the shape $\ll/\mu$ one step to the right. Similarly, define the
left shift of $\lambda/\mu$ by
$$
\la{\ll/\mu}:= (\ll_1-1 ,\ldots,\ll_k -1 )/(\mu_1-1,\ldots,\mu_k-1),
$$
assuming that the result is a legitimate skew shape.  Note that
$s_{\ll/\mu} = s_{\la{\ll/\mu}} = s_{\ra{\ll/\mu}}$.

Let $\theta$ be the operation on pairs of skew shapes given by
$$
\theta:(\lambda/\mu,\nu/\rho)\longmapsto ((\lambda/\mu)\vee(\nu/\rho),
(\lambda/\mu)\wedge(\nu/\rho)).
$$
According to Theorem~\ref{th:nonnegativity_s}, the product of the
two skew Schur functions corresponding to the shapes in
$\theta(\lambda/\mu,\nu/\rho)$ is $\geq_s
s_{\ll/\mu}\,s_{\nu/\rho}$. Let us show that we can repeatedly
apply the operation $\theta$ together with the left and right
shifts of shapes and the flips $(\lambda/\mu,\nu/\rho)\mapsto
(\nu/\rho,\lambda/\mu)$ in order to obtain the pair of skew shapes
$(\lfloor\frac{\ll+\nu}{2}\rfloor/\lfloor\frac{\mu+\rho}{2}\rfloor,
\lceil\frac{\ll+\nu}{2}\rceil/\lceil\frac{\mu+\rho}{2}\rceil)$
from $(\ll/\mu,\nu/\rho)$.

Let us define two operations $\phi$ and $\psi$
on ordered pairs of skew shapes by conjugating $\theta$ with
the right and left shifts and the flips, as follows:
\[
\begin{array}{l}
\phi:(\ll/\mu,\nu/\rho) \longmapsto ((\ll/\mu) \wedge (\ra{\nu/\rho}),
\la{(\ll/\mu) \vee (\ra{\nu/\rho})}),
\\[.1in]
\psi: (\ll/\mu,\nu/\rho) \longmapsto (\la{(\ra{\ll/\mu}) \vee (\nu/\rho)},
(\ra{\ll/\mu}) \wedge (\nu/\rho)).
\end{array}
\]
In this definition the application of the left shift ``$\la{}$'' always makes
sense.   Indeed, in both cases, before the application of ``$\la{}$'',
we apply ``$\ra{}$'' and then ``$\vee$''.
As we noted above, both products of skew
Schur functions for shapes in $\phi(\ll/\mu,\nu/\rho)$
and in $\psi(\ll/\mu,\nu/\rho)$ are $\geq_s s_{\ll/\mu}\,s_{\nu/\rho}$.

It is convenient to write the operations $\phi$ and $\psi$
in the coordinates $\lambda_i,\mu_i,\nu_i,\rho_i$, for $i=1,\dots,k$.
These operations independently act on the pairs $(\lambda_i,\nu_i)$ by
$$
\begin{array}{l}
\phi:(\lambda_i,\nu_i)\mapsto
(\min(\lambda_i,\nu_i+1),\max(\lambda_i,\nu_i+1)-1),\\
\psi:(\lambda_i,\nu_i)\mapsto
(\max(\lambda_i+1,\nu_i)-1,\min(\lambda_i+1,\nu_i)),\\
\end{array}
$$
and independently act on the pairs $(\mu_i,\rho_i)$
by exactly the same formulas.
Note that both operations $\phi$ and $\psi$ preserve the sums
$\lambda_i+\nu_i$ and $\mu_i+\rho_i$.

The operations $\phi$ and $\psi$ transform the differences
$\lambda_i-\nu_i$ and $\mu_i-\rho_i$ according to the following
piecewise-linear maps:
$$
\bar\phi(x) =\left\{
\begin{array}{cl}
x & \textrm{if } x\leq 1,\\
2-x & \textrm{if } x\geq 1,
\end{array}\right.
\qquad\textrm{and}\qquad
\bar\psi(x) =\left\{
\begin{array}{cl}
x & \textrm{if } x\geq -1,\\
-2-x & \textrm{if } x\leq -1.
\end{array}\right.
$$
Whenever we apply the composition $\phi\circ\psi$ of these operations,
all absolute values $|\lambda_i-\nu_i|$ and $|\mu_i-\rho_i|$ strictly
decrease, if these absolute values are $\geq 2$.
It follows that, for a sufficiently large integer $N$, we have
$(\phi\circ \psi)^N(\lambda/\mu,\nu/\rho) = (\tilde\lambda/\tilde\mu,
\tilde\nu/\tilde\rho)$ with
$\tilde\lambda_i+\tilde\nu_i = \lambda_i+\nu_i$,
$\tilde\mu_i+\tilde\rho_i = \mu_i+\rho_i$,
and
$|\tilde\lambda_i-\tilde\nu_i|\leq 1$, $|\tilde\mu_i-\tilde\rho_i|\leq 1$,
for all $i$.
Finally, applying the operation $\theta$, we obtain
$\theta(\tilde\lambda/\tilde\mu, \tilde\nu/\tilde\rho) =
(\lceil\frac{\ll+\nu}{2}\rceil/\lceil\frac{\mu+\rho}{2}\rceil,
\lfloor\frac{\ll+\nu}{2}\rfloor/\lfloor\frac{\mu+\rho}{2}\rfloor)$,
as needed.
\end{proof}

The following conjugate version of Theorem~\ref{th:twoshapes}
extends Fomin-Fulton-Li-Poon's conjecture (Conjecture~\ref{conj:FFLP})
to skew shapes.

\begin{corollary}
\label{cor:twoshapesconjugate}
Let $\ll/\mu$ and $\nu/\rho$ be
two skew shapes. Then we have
\[
s_{\mathrm{sort}_1(\ll,\nu)/\mathrm{sort}_1(\mu,\rho)}\,
s_{\mathrm{sort}_2(\ll,\nu)/\mathrm{sort}_2(\mu,\rho)} \geq_s
s_{\ll/\mu}\, s_{\nu/\rho}.
\]
\end{corollary}

\begin{proof}
This statement is obtained from Theorem~\ref{th:twoshapes} by conjugating
the shapes.
Indeed, $\lceil\frac{\lambda+\mu}{2}\rceil' = \mathrm{sort}_1(\lambda',\mu')$
and  $\lfloor\frac{\lambda+\mu}{2}\rfloor' =
\mathrm{sort}_2(\lambda',\mu')$.
Here $\lambda'$ denote the partition conjugate to $\lambda$.
\end{proof}

\begin{theorem}
\label{thm:skewlogconcavity2}
Let
$\ll^{(1)}/\mu^{(1)},\ldots,\ll^{(n)}/\mu^{(n)}$ be $n$ skew
shapes, let $\ll = \bigcup \ll^{(i)}$ be the partition obtained by
the decreasing rearrangement of the parts in all $\lambda^{(i)}$,
and, similarly, let $\mu = \bigcup\mu^{(i)}$.
Then we have
$
\prod_{i=1}^n s_{\ll^{[i,n]}/\mu^{[i,n]}} \geq_s \prod_{i=1}^n
s_{\ll^{(i)}/\mu^{(i)}}$.
\end{theorem}

This theorem extends Corollary~\ref{cor:twoshapesconjugate} and
Conjecture~\ref{conj:FFLP}.  Also note that
Lascoux-Leclerc-Thibon's conjecture (Conjecture~\ref{conj:LLT}) is
a special case of Theorem~\ref{thm:skewlogconcavity2} for the
$n$-tuple of partitions
$(\ll^{[1,m]},\ldots,\ll^{[m,m]},\emptyset,\ldots,\emptyset)$.

\begin{proof}
Let us derive the statement by applying
Corollary~\ref{cor:twoshapesconjugate}
repeatedly.
For a sequence  $v = (v_1,v_2,\ldots,v_l)$ of integers, the
{\it anti-inversion number} is $\ainv(v):=\#\{(i,j)\mid i<j,\ v_i<v_j\}$.
Let $L= (\ll^{(1)}/\mu^{(1)},\ldots,\ll^{(n)}/\mu^{(n)})$ be a sequence
of skew shapes.
Define its anti-inversion number as
$$
\begin{array}{l}
\ainv(L) = \ainv(\ll^{(1)}_1,\ll^{(2)}_1,\dots,\ll^{(n)}_1,\ll^{(1)}_2,
\dots,\ll^{(n)}_2, \ll^{(1)}_3,\dots,\ll^{(n)}_3,\dots)\\
\qquad\qquad
{}+\ainv(\mu^{(1)}_1,\mu^{(2)}_1,\dots,\mu^{(n)}_1,\mu^{(1)}_2,
\dots,\mu^{(n)}_2, \mu^{(1)}_3,\dots,\mu^{(n)}_3,\dots).
\end{array}
$$
If $\ainv(L) \ne 0$ then there is a pair $k<l$ such that
$\ainv(\ll^{(k)}/\mu^{(k)}, \ll^{(l)}/\mu^{(l)}) \ne 0$.
Let $\tilde L$ be the sequence of skew shapes obtained from $L$
by replacing the two terms $\ll^{(k)}/\mu^{(k)}$ and
$\ll^{(l)}/\mu^{(l)}$
with the terms
$$
\mathrm{sort}_1(\ll^{(k)}, \ll^{(l)})/
\mathrm{sort}_1(\mu^{(k)},\mu^{(l)})
\quad\textrm{and}\quad
\mathrm{sort}_2(\ll^{(k)}, \ll^{(l)})/
\mathrm{sort}_2(\mu^{(k)},\mu^{(l)}),
$$
correspondingly.  Then $\ainv(\tilde L)<\ainv(L)$.
Indeed, if we rearrange a subsequence in a sequence in
the decreasing order, the total number of anti-inversions decreases.
According to Corollary~\ref{cor:twoshapesconjugate},
we have $s_{\tilde L} \geq_s
s_L$, where $s_L := \prod_{i=1}^n s_{\ll^{(i)}/\mu^{(i)}}$.
Note that the operation $L\mapsto \tilde L$ does not change
the unions of partitions $\bigcup \lambda^{(i)}$ and $\bigcup \mu^{(i)}$.
Let us apply the operations $L\mapsto \tilde L$ for various pairs $(k,l)$ until
we obtain a sequence of skew shapes $\hat L=
(\hat\ll^{(1)}/\hat\mu^{(1)},\ldots,\hat\ll^{(n)}/\hat\mu^{(n)})$
with $\ainv(\hat L)=0$, i.e., the parts of all partitions must be sorted
as $\hat\lambda_1^{(1)}\geq \cdots \geq \hat\lambda_1^{(n)}\geq
\hat\lambda_2^{(1)}\geq \cdots \geq \hat\lambda_2^{(n)}\geq
\hat\lambda_3^{(1)}\geq \cdots \geq \hat\lambda_3^{(n)}\geq\cdots$,
and the same inequalities hold for the $\hat\mu_j^{(i)}$.
This means that $\hat\lambda^{(i)}/\hat\mu^{(i)} =
\lambda^{[i,n]}/\mu^{[i,n]}$, for $i=1,\dots,n$.
Thus $s_{\hat L} = \prod s_{\lambda^{[i,n]}/\mu^{[i,n]}}
\geq_s s_L$, as needed.
\end{proof}

Let us define $\lambda^{\{i,n\}} := ((\ll')^{[i,n]})'$, for
$i=1,\dots,n$. Here $\lambda'$ again denotes the partition
conjugate to $\lambda$. The partitions $\lambda^{\{i,n\}}$ are
uniquely defined by the conditions
$\lceil\frac{\lambda}n\rceil\supseteq \lambda^{\{1,n\}}
\supseteq\cdots\supseteq \lambda^{\{n,n\}} \supseteq
\lfloor\frac{\lambda}n\rfloor$ and $\sum_{i = 1}^n
\lambda^{\{i,n\}} = \ll$. In particular, $\lambda^{\{1,2\}}=\lceil
\frac{\lambda}2\rceil$ and $\lambda^{\{2,2\}}=\lfloor
\frac{\lambda}2\rfloor$. If $\frac{\lambda}n$ is a partition,
i.e., all parts of $\lambda$ are divisible by $n$, then
$\lambda^{\{i,n\}}=\frac{\lambda}n$ for each $1 \leq i \leq n$.

\begin{corollary}
\label{cor:lambda_plus_mu}
Let
$\ll^{(1)}/\mu^{(1)},\ldots,\ll^{(n)}/\mu^{(n)}$ be $n$ skew
shapes, let $\ll = \ll^{(1)}+\cdots+\lambda^{(n)}$
and $\mu = \mu^{(1)}+\cdots+\mu^{(n)}$.
Then we have
$\prod_{i=1}^n s_{\ll^{\{i,n\}}/\mu^{\{i,n\}}} \geq_s \prod_{i=1}^n
s_{\ll^{(i)}/\mu^{(i)}}$.
\end{corollary}

\begin{proof}
This claim is obtained from Theorem~\ref{thm:skewlogconcavity2}
by conjugating the shapes.  Indeed,
$\left(\bigcup \lambda^{(i)}\right)'
= \sum(\lambda^{(i)})'$.
\end{proof}

For a skew shape $\lambda/\mu$ and a positive integer $n$, define
$
s_{\frac{\ll}{n}/
\frac{\mu}{n}}^{\left<n\right>}:=
\prod_{i=1}^n s_{\ll^{\{i,n\}}/\mu^{\{i,n\}}}$.
In particular, if $\frac{\ll}{n}$ and $\frac{\mu}{n}$ are
partitions, then $s_{\frac{\ll}{n}/
\frac{\mu}{n}}^{\left<n\right>}= \left(s_{\frac{\ll}{n}/
\frac{\mu}{n}}\right)^n$.


\begin{corollary}
\label{cor:skewlogconcavity}
Let $c$ and $d$ be positive integers and $n=c+d$.
Let $\ll/\mu$ and $\nu/\rho$ be two skew shapes. Then
$
s_{\frac{c\ll + d\nu}{n}/\frac{c\mu + d\rho}{n}}^{\left<n\right>} \geq_s
s_{\ll/\mu}^c\, s_{\nu/\rho}^d$.
\end{corollary}

Theorem~\ref{th:twoshapes}
is a special case of Corollary~\ref{cor:skewlogconcavity} for
$c = d = 1$.

\begin{proof}
This claim follows from Corollary~\ref{cor:lambda_plus_mu} for the
sequence of skew shapes that consists of $\lambda/\mu$ repeated
$c$ times and $\nu/\rho$ repeated $d$ times.
\end{proof}

Corollary~\ref{cor:skewlogconcavity} implies that the map $S:\lambda
\mapsto s_\lambda$ from the set of partitions to symmetric
functions satisfies the following ``Schur log-concavity''
property.

\begin{corollary}
For positive integers $c,d$ and partitions
$\lambda,\mu$ such that $\frac{c\ll + d\mu}{c+d}$
is a partition, we have
$
\left(S\left(\frac{c\ll + d\mu}{c+d}\right)\right)^{c+d} \geq_s S(\ll)^c
S(\mu)^d.
$
\end{corollary}

This notion of Schur log-concavity is inspired by Okounkov's notion of
log-concavity; see~\cite{Oko}.

\bigskip
\noindent
{\sc Acknowledgements:}
We thank Richard Stanley for useful conversations.  
We are grateful to Sergey Fomin for helpful comments and suggestions
and to Mark Skandera for help with the references.

\end{document}